\documentclass[12pt]{amsart}
\usepackage{a4wide}
\usepackage{amssymb}
\usepackage{amsthm}
\usepackage{amsmath}
\usepackage{amscd}
\usepackage{verbatim}
\usepackage{stmaryrd}
\usepackage[all]{xy}
\usepackage{rotating}
\usepackage{rotfloat}

\numberwithin{equation}{section}

\theoremstyle{plain}
\newtheorem{theorem}{Theorem}[section]

\theoremstyle{definition}

\theoremstyle{remark}
\newtheorem{remark}{Remark}

\newcommand{\R}{\mathbb{R}}
\newcommand{\Q}{\mathbb{Q}}
\newcommand{\Z}{\mathbb{Z}}

\newcommand{\C}{\mathbb{C}}

\renewcommand{\H}{\mathbb{H}}

\newcommand{\leg}[2]{\left( \frac{#1}{#2} \right)}
\newcommand{\kzxz}[4]{\left(\begin{smallmatrix} #1 & #2 \\ #3 & #4\end{smallmatrix}\right) }
\newcommand{\kabcd}{\kzxz{a}{b}{c}{d}}

\newcommand{\frake}{\mathfrak e}

\newcommand{\eps}{\varepsilon}
\newcommand{\bs}{\backslash}

\newcommand{\Sl}{\operatorname{SL}}

\newcommand{\Mp}{\operatorname{Mp}}
\newcommand{\Orth}{\operatorname{O}}

\newcommand{\PSLZ}{\operatorname{PSL}_{2}(\mathbb{Z})}
\newcommand{\SLZ}{\operatorname{SL}_{2}(\mathbb{Z})}
\newcommand{\T}{T}

\makeatletter
\newskip\rotFPtop \rotFPtop=220pt
\newskip\rotFPbot \rotFPbot=0pt
\makeatother

\begin{document}

\title[Computation of harmonic weak Maass forms]{Computation of harmonic weak Maass forms}

\date{\today}
\author{Jan H. Bruinier and Fredrik Str\"omberg}

\address{Fachbereich Mathematik,
Technische Universit\"at Darmstadt, Schlossgartenstrasse 7, D--64289
Darmstadt, Germany} \email{bruinier@mathematik.tu-darmstadt.de}

\address{Fachbereich Mathematik,
Technische Universit\"at Darmstadt, Schlossgartenstrasse 7, D--64289
Darmstadt, Germany} \email{stroemberg@mathematik.tu-darmstadt.de}
\thanks{The authors are partially supported by DFG grant BR-2163/2-1. Hardware
obtained through the NSF grant DMS-0821725 has been used for a large part of the computations.}


\begin{abstract}
Harmonic weak Maass forms
of half-integral weight
are the subject of many recent works.
They are closely related to Ramanujan's mock theta functions, their theta lifts give rise to Arakelov Green functions, and their coefficients are often related to central values and derivatives of Hecke $L$-functions.
We present an algorithm to compute harmonic weak Maass forms numerically, based on the automorphy method due to Hejhal and Stark.
As explicit examples we consider harmonic weak Maass forms of weight $1/2$ associated to the elliptic curves 11a1, 37a1, 37b1.
We made extensive numerical computations and
the data we obtained is presented in the final section of the paper.
We expect that experiments based on our data will lead to a better understanding of the arithmetic properties of the Fourier coefficients.
\end{abstract}

\maketitle

\section{Introduction}
\label{sect:intro}

Half-integral weight modular forms play important roles in
arithmetic geometry and number theory. Their coefficients serve as generating functions for
various interesting number theoretic functions, such as representation numbers
of quadratic forms in an odd number of variables or class numbers of imaginary quadratic fields.
Moreover, employing the Shimura correspondence \cite{Sh},
Waldspurger \cite{Wa}, and Kohnen and Zagier
\cite{KZ,K} showed that
the coefficients of
half-integral weight cusp forms essentially are square-roots of
central values of quadratic twists of modular $L$-functions.
In analogy with these works, Katok and Sarnak \cite{KS} used a
Shimura correspondence to relate coefficients of weight 1/2 Maass
forms to sums of values and sums of line integrals of Maass cusp
forms.

In more recent work Zagier discovered that the generating function for the traces of
singular moduli (the CM values of the classical $j$-function) is a weakly holomorphic modular form of weight $3/2$ \cite{Za1}.
This result, which was generalized in various directions (see e.g. \cite{BO3}, \cite{BF2}, \cite{DJ}, \cite{Kim1}), demonstrates that also the coefficients of
automorphic forms with singularities at the cusps carry interesting arithmetic information.

In a similar spirit, Ono and the first author proved that the coefficients of harmonic weak Maass forms of weight $1/2$ are related to both the values
and central derivatives of quadratic twists of weight 2 modular
$L$-functions \cite{BruO}. Harmonic weak Maass forms are also closely related to mock modular forms and to Ramanujan's mock theta
functions, which have been the subject of various recent works (see e.g. \cite{BO1, BO2, O2, Za2, Z1, Z2}).
In view of these connections, it is desirable to develop tools for the computation of such automorphic forms.
In the present paper we propose an approach to this problem which yields an efficient algorithm.
Moreover, we compute some harmonic weak Maass forms which are related to rational elliptic curves as in \cite{BruO}.

The non-holomorphic nature of harmonic weak Maass forms prevents the use of the well developed algorithms existing for
(weakly) holomorphic modular forms, such as e.g. modular symbols.
The use of Poincar{\'e} series does not work well either in small weights due to the poor convergence of the infinite series which appear in the explicit formulas for the
coefficients.
Instead we adapt the  `automorphy method', originally developed by Hejhal for the computation of Maass cusp forms on Hecke triangle groups (see e.g. \cite{He}),
to the setting of harmonic weak Maass forms.

We now describe the content of this paper in more detail.
Let $k\in \frac{1}{2}\Z$, and let
$N$ be a positive integer (with $4\mid N$ if $k\in \frac{1}{2}\Z\setminus \Z$).
A {\em harmonic weak Maass form} of weight $k$ on
$\Gamma_0(N)$
is a smooth function on $\H$, the upper half of the complex plane,
which satisfies:
\begin{enumerate}
\item[(i)]
 $f\mid_k\gamma = f$ for all $\gamma\in \Gamma_0(N)$;
\item[(ii)] $\Delta_k f =0 $, where $\Delta_k$ is the weight $k$
hyperbolic Laplacian on $\H$ (see (\ref{deflap})); \item[(iii)]
There is a polynomial $P_f=\sum_{n\leq 0} c^+(n)q^n \in \C[q^{-1}]$
such that $f(\tau)-P_f(\tau) = O(e^{-\eps v})$ as $v\to\infty$ for
some $\eps>0$. Analogous conditions are required at all
cusps.
\end{enumerate}
Throughout, for $\tau\in \H$, we let $\tau=u+iv$, where $u, v\in
\R$, and we let $q:=e^{2 \pi i \tau}$.
The polynomial $P_f$ is called the {\em principal part} of $f$ at $\infty$.

Such a harmonic weak Maass form $f$ has a Fourier expansion at infinity of the form
\begin{equation}\label{fourier}
f(\tau)=\sum_{n\gg -\infty} c^+(n) q^n + \sum_{n<0} c^-(n)\Gamma\left(1-k, 4\pi |n|
v\right) q^n,
\end{equation}
where $\Gamma(a,x)$ denotes the incomplete Gamma function.
The series $\sum_{n\gg -\infty}
c^+(n) q^n$ is called the {\it holomorphic part} of $f$, and its
complement is called the {\it non-holomorphic part}. Naturally, $f$ has similar expansions at the other cusps.
There is an antilinear differential operator, taking $f$ to the cusp form
$\xi_{k}(f):=2iv^k\overline{\frac{\partial f}{\partial \bar \tau}}$ of weight $2-k$, see \eqref{defxi}.
The kernel of $\xi_k$ consists of the space of {\it weakly holomorphic}
modular forms, those meromorphic modular forms whose poles (if any)
are supported at cusps.

Every weight $2-k$ cusp form is the image under $\xi_{k}$ of a
weight $k$ harmonic weak Maass form. Ramanujan's mock theta functions
correspond  to those forms whose images under $\xi_{1/2}$
are weight 3/2 unary theta functions. Here we mainly consider
those weight 1/2 harmonic weak Maass forms whose images under
$\xi_{1/2}$ are orthogonal to the unary theta series.
According to \cite{BruO}, their coefficients
are related
to both the values
and central derivatives of quadratic twists of weight 2 modular
$L$-functions.

We now briefly describe this result in the special case that the level is a prime $p$.
Let
$G\in S_2(\Gamma_0(p))$ be
a normalized Hecke eigenform whose
Hecke $L$-function $L(G,s)$ satisfies an odd functional equation.
That is, the completed $L$-function $\Lambda(G,s)=p^{s/2}(2\pi)^{-s}\Gamma(s)L(G,s)$  satisfies $\Lambda(G,2-s)=\varepsilon_G \Lambda(G,s)$
with root number $\varepsilon_G=-1$.
Therefore, the central critical value $L(G,1)$ vanishes.
By Kohnen's theory of plus-spaces \cite{K}, there is a half-integral
weight newform
$
g\in S_{3/2}^{+}(\Gamma_0(4p))$,
unique up to a multiplicative constant, which lifts to $G$ under the
Shimura correspondence. We choose $g$ so that its
coefficients are in $F_G$, the
totally real number field generated by the Hecke eigenvalues  of $G$.
There exists a weight 1/2 harmonic weak Maass form $f$ on
$\Gamma_0(4p)$ in the plus space
whose principal part $P_{f}$ has coefficients in $F_G$, and such that
\[
\xi_{1/2}(f)={\|g\|^{-2}}g,
\]
where $\|g\|$ denotes the usual Petersson norm.
For a fundamental discriminant $\Delta$ let $\chi_{\Delta}$ be the
Kronecker character for $\Q(\sqrt{\Delta})$, and let $L(G,\chi_{\Delta},s)$ be
the quadratic twist of $L(G,s)$ by $\chi_{\Delta}$. One can show that the root number of $L(G,\chi_{\Delta},s)$ is equal to $\textrm{sign}(\Delta) \cdot \chi_{\Delta}(p)\, \varepsilon_G$.

\begin{theorem}[See \cite{BruO}]
\label{Lvalues}
Assume that $G$, $g$, and $f$ are as above, and let $c^{\pm}(n)$
denote the Fourier coefficients as in \eqref{fourier}.
\begin{enumerate}
 \item If $\Delta<0$ is a fundamental
discriminant for which $\leg{\Delta}{p}=1$, then
$$
L(G,\chi_{\Delta},1)=8\pi^2\|G\|^2 \|g\|^2 \sqrt{\frac{|\Delta|}{N}}\cdot
c^{-}(\Delta)^2.
$$
\item If $\Delta>0$ is a fundamental discriminant for which
$\leg{\Delta}{p}=1$, then $L'(G,\chi_{\Delta},1)=0$ if and only if
$c^{+}(\Delta)$ is algebraic.
\end{enumerate}
\end{theorem}

Note that the harmonic weak Maass form $f$ is uniquely determined up
to the addition of a weight 1/2 weakly holomorphic modular form with
coefficients in $F_G$.  Furthermore, the absolute values of the
nonvanishing coefficients $c^{+}(\Delta)$ are typically asymptotic to
subexponential functions in $n$. For these reasons, the connection
between $L'(G,\chi_{\Delta},1)$ and the coefficients $c^{+}(\Delta)$
in Theorem~\ref{Lvalues}(2) cannot be modified in a simple way to
obtain a formula as in the first part of the Theorem.  In fact, the
proof of Theorem~\ref{Lvalues}(2) is rather indirect. It relies on the
Gross-Zagier formula and on transcendence results of Waldschmidt and
Scholl on periods of differentials on algebraic
curves.

The above result is one of the main motivations for the present
paper. Our goal is to carry out numerical computations for the
involved harmonic weak Maass forms.  In that way we hope to find more
direct connections of the coefficients $c^+(\Delta)$ to periods or
$L$-functions.
When $L'(G,\chi_{\Delta},1)$ vanishes,
meaning that $c^+(\Delta)$ is algebraic (actually contained in $F_G$), it
would be interesting to see if $c^+(\Delta)$ carries any arithmetic
information related to $G$.
In a forthcoming paper \cite{Brperiods},
the coefficients $c^+(n)$ will be linked to periods of
certain algebraic differentials of the third kind on modular curves.
It leads to a conjecture on differentials of the third kind on elliptic curves,
which is based on the numerical data presented in Section \ref{sect:results} of the present paper.

Our computations make use of an adaption of the so-called automorphy
method. The key point of this method is to view an automorphic form on
a non-co-compact (but co-finite) Fuchsian group $\Gamma$ as a function
on the upper half-plane with certain transformation properties under
the group $\Gamma$ as well as convergent Fourier series expansions at
all cusps.  This classical point of view, in terms of functions on the
upper half-plane, stands in contrast to the more algebraic point of
view, in terms of Hecke modules, usually taken when computing
holomorphic modular forms.

By {\em computing} an automorphic form $\phi$ in this setting we mean
that to any given (small) $\epsilon > 0$ we compute a sufficient
number of Fourier coefficients, each to high enough precision, so that
we are able to evaluate the function $\phi$ at any point in the upper
half-plane with an error at most $\epsilon$.

To calculate these Fourier coefficients we truncate the Fourier series
representing $\phi$ and view the resulting trigonometric sum as a
finite Fourier series. Using the Fourier inversion theorem together
with the automorphic properties of $\phi$ (which will additionally
intertwine the Fourier series at various cusps) we are able to obtain
a set of linear equations satisfied approximately by the coefficients.
Cf.\,e.g.\,\cite{He,St,Av2}. The (surprising) effectiveness of this
algorithm is closely related to the equidistribution properties of
closed horocycles (cf.\, e.g.\,\cite{He1,S}).
We describe the main algorithm in detail in Section \ref{sect:comp}.
The implementation of the software package is briefly described in Section \ref{ssect:implementation}.

In Section \ref{sect:results} we describe our computational result in
three cases of particular interest.  We consider the elliptic curves
11a1, 37a1, and 37b1 and their corresponding weight $2$ newforms.  For
instance, the elliptic curve 37a1, is the curve of smallest conductor
with rank 1. It corresponds to the unique weight two normalized
newform $G$ on $\Gamma_0(37)$ whose $L$-function has an odd functional
equation.  We verified the statement of Theorem \ref{Lvalues} for all
fundamental discriminants $\Delta$ which are squares modulo $148$ in
the range $0<\Delta < 15000$.  For eight of these fundamental
discriminants the quantity $L'(G,\chi_{\Delta},1)$ vanishes.
In all these cases we
found a stronger statement then that of the Theorem \ref{Lvalues} to
be true, namely, that the associated coefficient
$c^+(\Delta)$ was an integer. For the corresponding data
see Tables \ref{tab:37a1} and \ref{tab:37a1-zeros}.
We conclude Section \ref{sect:results} by
describing some analogous experiments for newforms $G$ of weight $4$,
where $g$ is of weight $5/2$ and $f$ of weight $-1/2$.


The present paper is organized as follows.  In Section \ref{sect:2} we
recall some facts on (half integral weight) harmonic weak Maass
forms. When working with arbitrary (not necessarily prime) level, it
is convenient to use vector valued modular forms. In Section
\ref{sect:2.3} we therefore recall from \cite{BruO} the vector valued
version of Theorem \ref{Lvalues}.  In Section \ref{sect:comp} we
describe the automorphy method in the context of harmonic weak Maass
forms. In Section \ref{sect:results} we collect our computational
results. In particular, we present results for the elliptic curves
11a1, 37a1, and 37b1; cf., e.g.~Tables \ref{tab:11a1}, \ref{tab:37a1}
and \ref{tab:37b1}.  More extensive tables can be obtained from the
authors on request.


\section{Preliminaries}
\label{sect:2}

In order to be able to work with newforms of arbitrary level, it is convenient to work with vector valued modular forms of half integral weight for the metaplectic extension of $\Sl_2(\Z)$. We describe the necessary background in this section.

\subsection{A Weil representation}

Let $\H=\{\tau\in \C;\;\Im(\tau)>0\}$ be the complex upper half
plane. We write $\Mp_2(\R)$ for the metaplectic two-fold cover of
$\Sl_2(\R)$, realized as the group of pairs $(M,\phi(\tau))$,
where $M=\kabcd\in\Sl_2(\R)$ and $\phi:\H\to \C$ is a holomorphic
function with $\phi(\tau)^2=c\tau+d$.  The multiplication is defined
by
\[
(M,\phi(\tau)) (M',\phi'(\tau))=(M M',\phi(M'\tau)\phi'(\tau)).
\]
We denote  the inverse image of
$\Gamma:=\Sl_2(\Z)$ under the covering map by
$\tilde\Gamma:=\Mp_2(\Z)$. It is well known that $\tilde\Gamma$ is
generated by $T:= \left( \kzxz{1}{1}{0}{1}, 1\right)$, and $S:=
\left( \kzxz{0}{-1}{1}{0}, \sqrt{\tau}\right)$.

Let $N$ be a positive integer. There is a certain representation $\rho$ of $\tilde \Gamma$ on
$\C[\Z/2N\Z]$, the group ring of the finite cyclic group of order $2N$.
For a coset $h\in \Z/2N\Z$ we denote by $\frake_h$ the corresponding standard basis vector of $\C[\Z/2N\Z]$. We write $\langle\cdot,\cdot \rangle$ for
the standard scalar product (antilinear in the second entry) such
that $\langle \frake_h,\frake_{h'}\rangle =\delta_{h,h'}$.
In terms of the generators $T$ and $S$ of $\tilde \Gamma$, the representation $\rho$ is given by
\begin{align}
\label{eq:weilt}
\rho(T)(\frake_h)&=e\left(\frac{h^2}{4N}\right)\frake_h,\\
\label{eq:weils}
\rho(S)(\frake_h)&=
\frac{1}{\sqrt{2iN}} \sum_{h' \; (2N)} e\left(-\frac{hh'}{2N}\right)
 \frake_{h'}.
\end{align}
Here the sum runs through the elements of $\Z/2N\Z$ and we have put $e(a)=e^{2\pi i a}$.
Note that $\rho$ is the Weil representation associated to the one-dimensional positive definite lattice $K=(\Z,Nx^2)$ in the sense of \cite{Bor1}, \cite{Br}, \cite{BruO}.
It is unitary with respect to the standard scalar product.

If $k\in \frac{1}{2}\Z$, we write $M^!_{k,\rho}$ for the space of $\C[\Z/2N\Z]$-valued weakly holomorphic modular forms of weight $k$ for $\tilde \Gamma$ with representation $\rho$.
The subspaces of holomorphic modular forms and cusp forms are denoted by $M_{k,\rho}$ and $S_{k,\rho}$, respectively.

\subsection{Harmonic weak Maass forms}

\label{sect:2.2}
In this subsection
we assume that
$k\leq 1$.  A twice continuously differentiable function $f:\H\to
\C[\Z/2N\Z]$ is called a {\em harmonic weak Maass form} (of weight $k$ with
respect to $\tilde \Gamma$ and $\rho$) if it satisfies:
\begin{enumerate}
\item[(i)]
$f(M\tau) = \phi(\tau)^{2k}\rho(M,\phi) f(\tau)$
for all $(M,\phi)\in \tilde\Gamma$;
\item[(ii)]
$\Delta_k f=0$,
\item[(iii)]
there is a $\C[\Z/2N\Z]$-valued Fourier polynomial
\[
P_f(\tau)=\sum_{h\;(2N)}\sum_{n\in \Z_{\leq 0}} c^+(n,h) q^{\frac{n}{4N}} \frake_h
\]
such that $f(\tau)-P_f(\tau)=O(e^{-\eps v})$ as $v\to \infty$ for
some $\eps>0$.
\end{enumerate}
%
Here we have that
\begin{equation}
\label{deflap}
\Delta_k := -v^2\left( \frac{\partial^2}{\partial u^2}+
\frac{\partial^2}{\partial v^2}\right) + ikv\left(
\frac{\partial}{\partial u}+i \frac{\partial}{\partial v}\right)
\end{equation}
is the usual weight $k$ hyperbolic Laplace operator (see \cite{BF}).
The Fourier polynomial $P_f$  is called the {\em principal part} of
$f$. We denote the vector space of these harmonic weak Maass forms
by  $H_{k,\rho}$ (it was called $H^+_{k,\rho}$ in \cite{BF}).
Any weakly holomorphic modular form is a harmonic weak Maass form.
The Fourier expansion of any $f\in H_{k,\rho}$ gives a unique
decomposition $f=f^++f^-$, where
\begin{subequations}
\label{deff}
\begin{align}
\label{deff+}
f^+(\tau)&= \sum_{h\;(2N)}\sum_{\substack{n\in \Z\\ n\gg-\infty}} c^+(n,h) q^{\frac {n}{4N}}\frake_h,\\
\label{deff-}
f^-(\tau)&= \sum_{h\; (2N)}\sum_{\substack{n\in \Z\\ n< 0}} c^-(n,h) \Gamma\left(1-k,4\pi\left|\frac{n}{4N}\right|v\right) q^{\frac{n}{4N}} \frake_h.
\end{align}
\end{subequations}
We refer to $f^+$ as the {\em holomorphic part} and to $f^-$ as the {\em non-holomorphic part} of
$f$. Note that $c^{\pm}(n,h)=0$ unless $n\equiv h^2\, (4N)$.

Recall that there is an antilinear  differential operator $\xi=
\xi_k:H_{k,\rho}\to S_{2-k,\bar\rho}$, defined by
\begin{equation}
\label{defxi} f(\tau)\mapsto \xi(f)(\tau):=2iv^k\overline{\frac{\partial f}{\partial \bar \tau}}.
\end{equation}
Here $\bar\rho$ denotes the complex conjugate of the representation
$\rho$, which can be identified with the dual representation.  The map
$\xi$ is surjective and its kernel is the space $M^!_{k,\rho}$.
There is a bilinear pairing between
$M_{2-k,\bar\rho}$ and $H_{k,\rho}$ defined by the Petersson scalar
product
\begin{equation}\label{defpair}
\{g,f\}=\big( g,\, \xi(f)\big)
:=\int_{\Gamma\bs \H}\langle g,\, \xi(f)\rangle v^{2-k}\frac{du\,dv}{v^2},
\end{equation}
for $g\in M_{2-k,\bar\rho}$ and $f\in H_{k,\rho}$. If $g$ has
the Fourier expansion $g=\sum_{h,n} b(n,h) q^{n/4N}\frake_h$, and if we
denote the Fourier expansion of $f$  as in \eqref{deff}, then by
 \cite[Proposition 3.5]{BF} we have
\begin{equation}\label{pairalt}
\{g,f\}= \sum_{h\;(2N)} \sum_{n\leq 0}  c^+(n,h) b(-n,h).
\end{equation}
Hence $\{g,f\}$ only depends on the principal part of $f$.

\subsection{The Shimura lift}
\label{sect:2.3}

Let $k\in
\frac{1}{2}\Z\setminus \Z$.  According to \cite[Chapter 5]{EZ}, the
space $M_{k,\bar\rho}$ is isomorphic to $J_{k+1/2,N}$, the space of
holomorphic Jacobi forms of weight $k+1/2$ and index $N$. According to
\cite{Sk1} and \cite{SZ}, $M_{k,\rho}$ is isomorphic to
$J_{k+1/2,N}^{skew}$, the space of skew holomorphic Jacobi forms of
weight $k+1/2$ and index $N$.  There is an extensive Hecke theory for
Jacobi forms (see \cite{EZ}, \cite{Sk1}, \cite{SZ}), which gives rise
to a Hecke theory on $M_{k,\rho}$ and $M_{k,\bar\rho}$, and which
is compatible with the Hecke theory on vector valued modular forms
considered in \cite{BrSt}.  In particular, there is an Atkin-Lehner
theory for these spaces.

The subspace $S_{k,\rho}^{new}$ of newforms of $S_{k,\rho}$ is
isomorphic as a module over the Hecke algebra to the space of
newforms $S^{new,+}_{2k-1}(N)$ of weight $2k-1$ for $\Gamma_0(N)$ on
which the Fricke involution acts by multiplication with
$(-1)^{k-1/2}$. The isomorphism is given by the Shimura
correspondence. Similarly, the subspace $S_{k,\bar\rho}^{new}$ of
newforms of $S_{k,\bar\rho}$ is isomorphic as a module over the
Hecke algebra to the space of newforms $S^{new,-}_{2k-1}(N)$ of
weight $2k-1$ for $\Gamma_0(N)$ on which the Fricke involution acts
by multiplication with $(-1)^{k+1/2}$ (see \cite{SZ}, \cite{GKZ},
\cite{Sk1}). Observe that the Hecke $L$-series of any $G\in
S^{new,\pm}_{2k-1}(N)$ satisfies a functional equation under
$s\mapsto 2k-1-s$ with root number $\varepsilon_G=\pm 1$.

We now state the vector valued version of Theorem \ref{Lvalues}.
Let $G\in
S_{2}^{new}(N)$ be a normalized newform (in particular a
common eigenform of all Hecke operators) of weight $2$ and write
$F_G$ for the number field generated by the eigenvalues of $G$. If
$\eps_G=-1$ we put $\rho'=\rho$, and
if $\eps_G=+1$ we put $\rho'=\bar \rho$.
There is a newform
$g\in S_{3/2,\bar\rho'}^{new}$ mapping to $G$ under the Shimura
correspondence. It is well
known that we may normalize $g$ such that all its coefficients
are contained in $F_G$. According to \cite[Lemma 7.3]{BruO}, there is a
harmonic weak Maass form $f\in H_{1/2, \rho'}$ whose principal part has coefficients in $F_G$ with the property that
\[
\xi_{1/2}(f)=\|g\|^{-2} g.
\]
This form is  unique up to addition of a weakly holomorphic form in
$M^!_{1/2, \rho'}$ whose principal part has coefficients in $F_G$.

In practice, the principal part of such an $f$ can be computed as follows:
We may complete the weight $3/2$ form $g$ to an orthogonal basis $g, g_2,\dots , g_d$ of $S_{3/2,\bar\rho'}$ consisting of cusp forms with Fourier coefficients in $F_G$. Let $f\in H_{1/2, \rho'}$ such that
\begin{align}
\label{eq:conds}
\{f,g\}=1, \quad \text{and $\{f,g_i\}=0$ for $i=2,\dots d$}.
\end{align}
Then $f$ has the required properties.
In view of \eqref{pairalt} the conditions of \eqref{eq:conds} translate into an inhomogeneous system of linear equations for the principal part of $f$.

\begin{theorem}
\label{Lval2}
Let $G\in S_{2}^{new}(N)$ be a normalized newform.
 Let $g\in S_{3/2,\bar\rho'}^{new}$, and $f\in H_{1/2, \rho'}$ be as above. Denote the Fourier coefficients of $f$ by $c^\pm(n,h)$ for $n\in \Z$ and $h\in \Z/2N\Z$. Then the following are true:
\begin{enumerate}
\item If $\Delta\neq 1$ is a fundamental
discriminant and $r\in \Z$ such that $\Delta\equiv r^2\pmod{4N}$ and $\varepsilon_G \Delta >0$, then
$$
L(G,\chi_{\Delta},1)=8\pi^2\|G\|^2 \|g\|^2 \sqrt{\frac{|\Delta|}{N}}\cdot
c^{-}(\Delta)^2.
$$
\item If $\Delta\neq 1$ is a fundamental
discriminant and $r\in \Z$ such that $\Delta\equiv r^2\pmod{4N}$ and $\varepsilon_G \Delta <0$, then
$$
L'(G,\chi_{\Delta},1)=0 \quad \Longleftrightarrow\quad c^{+}(-\varepsilon_G\Delta,r)\in \bar \Q\quad   \Longleftrightarrow\quad c^{+}(-\varepsilon_G\Delta,r)\in F_G.$$
\end{enumerate}
\end{theorem}

When $S_{1/2,\rho'}=\{0\}$ the above result also holds for $\Delta=1$, see also \cite[Remark 18]{BruO}. This is for instance the case when $N$ is a prime.
If $N$ is a prime and $\eps_G=-1$, then the space $H_{1/2, \rho'}$ can be identified with a space of scalar valued modular forms satisfying a Kohnen plus space condition.
In that way one obtains Theorem \ref{Lvalues} stated in the introduction.

\section{Computational aspects}
\label{sect:comp}

\subsection{The automorphy method for vector valued weak Maass forms}

To compute the Fourier coefficients of the harmonic weak Maass forms
we use the so-called automorphy method, sometimes called ``Hejhal's
method{}''. This is a general method which has been used to successfully
compute various kinds of automorphic functions and forms on $\mbox{GL}_{2}\left(\mathbb{R}\right)$.
It was originally developed by Hejhal in order to compute Maass cusp forms
for the modular group and other Hecke triangle groups (cf.~e.g.~\cite{He}).
The method was later generalized by the second author in \cite{St}
to computations of Maass waveforms with non-trivial multiplier systems
and arbitrary real weights, as well as to general subgroups of the
modular group (see also \cite{St2}).
Another generalization to automorphic forms with singularities
(Eisenstein series, Poincar{\'e} series and Green's functions) was made
by Avelin \cite{Av1,Av2}.

We will detail the adaptation of the algorithm to the case of vector-valued
harmonic weak Maass forms for the Weil representation.

For simplicity consider the representation $\rho$ (the case of $\overline{\rho}$
is analogous) and $k\in\mathbb{Z}+\frac{1}{2}$. Furthermore, in order
to avoid questions of uniqueness we assume that either $k<0$ or that $k=\frac{1}{2}$ and that $N$ is prime.
In these cases, a harmonic weak Maass form is uniquely determined by its principal part.
For computational purposes it is not feasible to use the definition of $\rho$ in terms
of the action on the generators of the metaplectic group. We instead
use formulas from {[}St1{]} to evaluate $\rho$ on the fixed (canonical)
representative of $M=\left(\begin{smallmatrix}a & b\\
c & d\end{smallmatrix}\right)\in\SLZ$, i.e. $\rho\left(M\right):=\rho\left(M,j_{M}\left(\tau\right)\right)$
where $j_{M}\left(\tau\right)=\sqrt{c\tau+d}$ is defined by the principal
branch of the argument.

\subsubsection{The algorithm -- phase 1}

Let $f\in H_{k,\rho}$ with a given (fixed) principal part $P_{f}\left(\tau\right)=\sum_{h}P_{f,h}\left(\tau\right)\mathfrak{e}_{h}$
where $P_{f,h}\left(\tau\right)=\sum_{n=-K}^{0}a\left(n,h\right)q^{\frac{n}{4N}}$
(for some finite $K\ge0$) and write $f=f^{+}+f^{-}$ (as in 2.3a and 2.3b)
with $f^{+}=\sum_{h\left(2N\right)}f_{h}^{+}\mathfrak{e}_{h}$ and
$f^{-}=\sum_{h\left(2N\right)}f_{h}^{-}\mathfrak{e}_{h}$ where
\begin{align*}
f_{h}^{+}\left(\tau\right) & =  \sum_{n=-K}^{0}a\left(n,h\right)q^{\frac{n}{4N}}+\sum_{n>0}c^{+}\left(n,h\right)q^{\frac{n}{4N}}\quad\mbox{and}\\
f_{h}^{-}\left(\tau\right) & =  \sum_{n<0}c^{-}\left(n,h\right)\Gamma\left(1-k,4\pi\left|\frac{n}{4N}\right|v\right)q^{\frac{n}{4N}}
\end{align*}
for $\tau=u+iv\in\H$. Our goal is to obtain numerical approximations to the coefficients $c^{\pm}(n,h)$.
To formulate our algorithm we prefer to separate the $u$- and the $v$-dependence in $f$ and therefore introduce the function $W$ defined by
$W(v)=e^{-2\pi v}$ if $v>0$ and $W(v)=e^{-2\pi v}\Gamma(1-k,4\pi|v|)$ if $v<0$.
We also set $c\left(n,h\right)=c^{+}\left(n,h\right)$ for ~$n>0$ and $c^{-}\left(n,h\right)$ for $n<0$ and write
$e_{4N}\left(u\right)=e^{\frac{2\pi iu}{4N}}$. With this notation
\[
f_{h}\left(\tau\right)=\sum_{n=-K}^{0}a\left(n,h\right) q^{\frac{n}{4N}}+
\sum_{n\ne0}c\left(n,h\right) W\left(\frac{nv}{4N}\right) e_{4N}(nu).
\]
By standard inequalities for the incomplete gamma function one can
show that \[
|W(v)|<c_{k}\, e^{-2\pi|v|}
\begin{cases}
1, & v>0,\\
\left(4\pi\left|v\right|\right)^{-k}, & v<0,
\end{cases}\]
where $c_{k}$ is an explicit constant only depending on $k$. To
be able to determine a truncation point of the Fourier series above
we also need bounds of the coefficients $c\left(n,h\right)$.
Using {[}BruFu, Lemma 3.4{]} it follows that there exists an explicit constant $C>0$ such that
\begin{align*}
c\left(n,h\right) & =  O\left(\exp\left(4\pi C\sqrt{n}\right)\right),\quad n\rightarrow+\infty,\\
c\left(n,h\right) & =  O(|n|^{\frac{k}{2}}),\quad n\rightarrow-\infty.
\end{align*}
For $k<0$ we are able to make the implied constants explicit using non-holomorphic Poincar{\'e}
series as in e.g.~\cite{Br} or \cite{He2}. For $k=\frac{1}{2}$ we rely on numerical a posteriori tests to assure ourselves that the truncation point was choosen correctly.
See e.g.~Section \ref{heuristics}.

Let $\epsilon>0$
and fix $Y<Y_{0}=\frac{\sqrt{3}}{2}$. By the estimates above we can
find an $M_{0}=M\left(Y,\epsilon \right)$ such that the function $\hat{f}=\sum_{h\,(2N)} \hat{f}_{h} \mathfrak{e}_h$ given by the truncated Fourier series
\[
\hat{f}_{h}\left(\tau\right)=P_{f,h}\left(\tau\right)+\sum_{0<\left|n\right|\le M_{0}}c\left(n,h\right)W\left(\frac{nv}{4N}\right)e_{4N}\left(nu\right)
\]
satisfies
\[
 \left\Vert \hat{f}(\tau) - f(\tau) \right\Vert^{2} < \epsilon
\]
for any $\tau \in \mathcal{H}_{Y}=\left\{\tau \in \mathcal{H}\,|\,\Im \tau \ge Y \right\}$. Here $\left\Vert z \right\Vert^{2}=\sum_{h=1}^{2N} \left| z_h \right|^2$ for $z \in \mathbb{C}^{2N}$.
Let $A=\left(\begin{smallmatrix} a & b \\ c & d \end{smallmatrix}\right)\in \SLZ$  and set $z=x+iy = A \tau$. Then
$
 y = \Im A\tau = \frac{v}{\left|c\tau + d\right|^2} \le \frac{v}{c^2v^2}\le \frac{1}{v}
$
and hence $\left|j_A(\tau)\right|^{4} = \left|c\tau + d\right|^{2} = \frac{v}{y} \le \frac{1}{y^2}$.
Using the fact that $\rho$ is unitary it is now easy to see that if  $\tau, A\tau \in \mathcal{H}_{Y}$ then
\begin{align}
\label{approx-autom}
\left\Vert \hat{f}\left(A\tau\right) -  j_{A}\left(\tau\right)^{2k}\rho\left(A\right)\hat{f}\left(\tau\right) \right\Vert^{2} < \epsilon \left( 1 + Y^{-2k} \right) < 2 \epsilon \cdot Y^{-2k}.
\end{align}
Consider now a horocycle at height $Y$ and a set of $2Q$ (with $Q>M_{0}$) equally spaced
points
\[
z_{m}=x_{m}+iY,\quad x_{m}=\frac{1-2m}{4Q},\quad1-Q\le m\le Q.
\]
If we view the series $\hat{f}_{h}$ as a finite Fourier series we can invert it over this horocycle and it is easy to see
that if $n$ is an integer with $0<\left|n\right|\le M_{0}$ and $n\equiv h^{2} \,(4N)$ then
\begin{equation}
\frac{1}{2Q}\sum_{m=1-Q}^{Q}\hat{f}_{h}\left(z_{m}\right)e_{4N}\left(-nx_{m}\right)=W\left(\frac{n}{4N} Y\right)c\left(n,h\right)+a\left(n,h\right)\,e^{-\frac{2\pi n}{4N}Y}.\label{eq:ffseries}
\end{equation}
One can also interpret the left-hand side as a Riemann-sum approximation
to the integral
\[\int_{-\frac{1}{2}}^{\frac{1}{2}}f_{h}\left(z\right)e_{4N}\left(-nx\right)dx.
\]
Let $z_{m}^{*}=x_{m}^{*}+iy_{m}^{*}=\T_{m}^{-1}z_{m}$ ($\T_{m}\in\PSLZ$)
denote the pull-back of $z_{m}$ into the standard (closed) fundamental domain of $\PSLZ$,
 $\mathcal{F}=\left\{ z=x+iy\,|\,\left|x\right|\le\frac{1}{2},\,\left|z\right|\ge1\right\} $.
Using (\ref{approx-autom}) we obtain
\[
\hat{f}_{h}\left(z_{m}\right)  = j_{\T_{m}}(z_{m}^{*}) \sum_{h'\,(2N)}\rho_{hh'}\left(\T_{m}\right)\hat{f}_{h'}\left(z_{m}^{*}\right) + \llbracket 2\epsilon Y^{-2k} \rrbracket,
\]
where $\rho_{hh'}\left(T_{m}\right)$ is the $\left(h,h'\right)$-element
of the matrix $\rho\left(T_{m}\right)$, and we use $\llbracket 2\epsilon Y^{-2k} \rrbracket$ to denote a quantity bounded in absolute value by $2\epsilon Y^{-2k}$.
Inserting this into (\ref{eq:ffseries}) we see that the left-hand side can be written as
\begin{align}
 \frac{1}{2Q} \sum_{m=1-Q}^{Q} j_{ \T_{m}} \left(z_{m}^{*}\right)  \sum_{h'\,(2N)} \rho_{hh'}\left(\T_{m}\right)
 & \left[ \sum_{l=-K}^{0} a (l,h') \exp\left(-\frac{2\pi l}{4N} y_{m}^{*}\right) e_{4N}\left(l x_{m}^{*}\right) \right. \nonumber \\
  + & \left.\sum_{0 < |l| \le M_{0} } c(l,h') W \left(\frac{l}{4N} y_{m}^{*}\right) e_{4N} \left(l x_{m}^{*}\right) \right] e_{4N}(-n x_{m})\nonumber \\
  = & \sum_{h'\,(2N)}\sum_{0<\left|l\right|\le M_{0}} c\left(l,h'\right) \widetilde{V}_{nl}^{hh'}+\widetilde{W}_{n}^{h}  + \llbracket 2\epsilon Y^{-2k} \rrbracket,
\label{eq:fundamental_expression}
\end{align}
where \begin{align*}
\widetilde{V}_{nl}^{hh'} & =  \frac{1}{2Q}\sum_{m=1-Q}^{Q}j_{T_{m}}\left(z_{m}^{*}\right)\rho_{hh'}\left(T_{m}\right)W\left(\frac{l}{4N} y_{m}^{*}\right)e_{4N}(lx_{m}^{*}-nx_{m})\quad\mbox{and}\\
\widetilde{W}_{n}^{h} & =  \frac{1}{2Q}\sum_{h'\,(2N)}\sum_{l=-K}^{0}a\left(l,h'\right)\sum_{m=1-Q}^{Q}j_{T_{m}}\left(z_{m}^{*}\right) \rho_{hh'}\left(T_{m}\right)
\exp\left(-\frac{2\pi l}{4N} y_{m}^{*}\right) e_{4N}\left(lx_{m}^{*}-nx_{m}\right).\end{align*}
We thus have an inhomogeneous system of linear equations which is (approximately)
 satisfied by the coefficients $c\left(n,h\right)$. Let $\mathcal{D}=\left\{ \left(n,h\right)\,|\,0<\left|n\right|\le M_{0},\,0\le  h < 2N \right\} $
(with a fixed ordering) and note that $\left|\mathcal{D}\right|=4M_{0}N$.
If we set $\vec{D}=\left(d\left(n,h\right)\right)_{\left(n,h\right)\in\mathcal{D}}$,
\begin{align*}
V &= V\left(Y\right)=\left(V_{nl}^{hh'}\right)_{\left(h,n\right),\left(h',l\right)\in\mathcal{D}},&
V_{nl}^{hh'}&=\widetilde{V}_{nl}^{hh'}-\delta_{nl}\delta_{hh'}W\left(\frac{n}{4N}Y\right)  \quad\text{and} \\
\vec{W} &= \vec{W}\left(Y\right)=\left(W_{n}^{h}\right)_{\left(h,n\right)\in\mathcal{D}}, &
W_{n}^{h}&=\widetilde{W}_{n}^{h}-a\left(n,h\right)e^{-\frac{2\pi n}{4N}Y},
\end{align*}
we can write this linear system as $\left|\mathcal{D}\right|$ linear equations in $\left|\mathcal{D}\right|$ variables:
\begin{equation}
V\vec{D}+\vec{W}=\vec{0}.\label{eq:linsyst}\end{equation}
In practice it turns out that the the matrix $V$ is non-singular as soon as the subspace of $H_{k,\rho}$ consisting of functions with a given singular part is one-dimensional.
In these cases we can immediately obtain the solution as
 \[
\vec{D}=-V^{-1}\vec{W},
\]
and since we know that the vector of the ``true'' coefficients, $\vec{C}=\left(c\left(n,h\right)\right)_{\left(n,h\right)\in\mathcal{D}}$, satisfies
\[
 \left\Vert V\vec{C}+\vec{W} \right\Vert_{\infty} \le 2\epsilon Y^{-2k},
\]
we see that
\[
 \left\Vert \vec{C}-\vec{D} \right\Vert_{\infty} =  \left\Vert \vec{C}+V^{-1} \vec{W} \right\Vert_{\infty} \le
 \left\Vert V^{-1} \right \Vert_{\infty} \cdot \left\Vert V\vec{C}+\vec{W} \right\Vert_{\infty} \le 2\epsilon Y^{-2k} \left\Vert V^{-1} \right \Vert_{\infty}.
\]
To obtain a theoretical error estimate we would thus need to estimate $\left\Vert V^{-1} \right \Vert_{\infty}$ from below.
Unfortunately this does not seem to be possible from the formulas above and we have to use numerical methods to estimate this norm.
Hence, to obtain the Fourier coefficients up to a (proven) desired precision we might have to go back and decrease the original $\epsilon$ or increase either of $M_0$ or $Q$.

At this point one should also remark that the error bound  $\left\Vert V^{-1} \right \Vert_{\infty}$ is in general much worse than the actual apparent error,
as verified by studying coefficients known to be integers. The reason for this is that the sums $\widetilde{V}_{nl}^{hh'}$ exhibit  massive cancellation
and are therefore overpowered by the terms $W\left( \frac{n}{4N} Y\right)$ on the diagonal.

\subsubsection{The algorithm -- phase 2}
Returning to (\ref{eq:fundamental_expression}) and solving for $c\left(n,h\right)$
we see that
\begin{equation}
c\left(n,h\right)=W\left(\frac{n}{4N} Y\right)^{-1}\left[\sum_{h'\,(2N)}\sum_{\left|l\right|\le M_{0}}c\left(l,h'\right)\widetilde{V}_{nl}^{hh'}+W_{n}^{h}   + \llbracket 2\epsilon Y^{-2k} \rrbracket \right]\label{eq:phase2}
\end{equation}
for \emph{any} $n$, i.e. also when $\left|n\right|>M_{0}$, provided
that $Q>M\left(Y\right)$. If we first choose $Y$ such that $W\left(\frac{n}{4N} Y\right)$
is not too small then we can in fact use this equation to compute
$c\left(n,h\right)$ with an error of size $\epsilon\,W\left(\frac{n}{4N} Y\right)^{-1}$.
In this manner, we may produce long stretches of coefficients (before
we need to decrease $Y$ again) at arbitrary intervals $N_{A}\le n\le N_{B}$
without the need of computing intermediate coefficients above the
initial set up to $n=M_{0}$.

\begin{remark}
The exact same algorithm, with the non-holomorphic parts set to zero, also lets one
compute holomorphic vector-valued modular forms for the Weil representation.
This has been exploited by the second author, in verifying computations of holomorphic Poincar\'e series in \cite{RSS}.
\end{remark}

\subsection{Heuristic error estimates}
\label{heuristics}

For $k<0$ all implied constants and therefore all error estimates can be made explicit.
In the remaining case which interests us, $k=\frac{1}{2}$, the known bounds for the twisted Kloosterman sums are not enough to
prove the necessary explicit bounds for the Fourier coefficients of the associated Poincar\'e series. We are therefore not able to give
effective theoretical error estimates in this case.
However, this is not a serious problem since there are a number of tests we may perform on the resulting coefficients
to assure ourselves of their accuracy. We list a few tests which we have used.
\begin{itemize}
 \item First of all, one can simply use two different values of $Y$ and verify that the resulting vectors $\vec{D}=\vec{D}(Y)$
are independent of $Y$.
\end{itemize}
This test is completely general and can be used for all instances where the algorithm can be applied.
Suppose now that we have a harmonic weak Maass form $f\in H_{k,\rho}$ of half-integral weight $k$ such that
$\xi_{k}\left(f\right)=\|g\|^{-2} g$, with $g\in S_{2-k,\bar\rho}$.
We then know the following.
\begin{itemize}
 \item The coefficients $\sqrt{\left|\Delta\right|} c^{-}(-\varepsilon_G\cdot\Delta)$ are proportional
to the coefficients $b(\varepsilon_G\cdot\Delta)$ of $g$ (cf.~e.g.~\cite[p.~3]{BruO}).
\end{itemize}
If additionally the Shimura lift of $g$ is a newform $G\in S^{new}_{3-2k}\left(\Gamma_0(N)\right)$
then we can predict that certain coefficients $c^{+}(\Delta)$ are algebraic (cf.~e.g.~\cite[Sect.~7]{BruO}) and if
we are able to identify these coefficients as algebraic numbers to a certain precision
this can be used as another measure of the accuracy.

\subsection{Implementation}
\label{ssect:implementation}
The first implementation of the above described algorithm was made in Fortran 90, using the package ARPREC \cite{AR} for arbitrary (fixed) precision arithmetic.
The second and more recent implementation was done in Sage \cite{SA}, using the included package mpmath for arbitrary (fixed) precision arithmetic.
The algorithms are currently under development but can be obtained on request from the authors.
The final format we intend for these algorithms are standard classes for computing with vector and scalar-valued harmonic weak Maass forms in Sage or Purple Sage.

\section{Results}
\label{sect:results}

\subsection{Harmonic Maass forms corresponding to elliptic curves}
In this section we present the numerical results we have obtained for harmonic weak Maass forms corresponding to weight two holomorphic forms associated to elliptic curves.
We have concentrated on three particular examples. In Cremona's notation, these correspond to the curve $11a1$ of level 11 and the two curves $37a1$ and $37b1$ of level 37.

Recall that if the holomorphic weight $2$ newform $G$ of level $N$ has Atkin-Lehner eigenvalue $\pm1$ then the $L$-function $L(G,s)$ has root number $\varepsilon_G=\mp1$.
Furthermore, since the root number of the twisted $L$-function $L(G,\chi_{\Delta},s)$ is $\textrm{sign}(\Delta)\chi_{\Delta}(N)\varepsilon_G$ and we always consider fundamental  discriminants for which $\chi_\Delta(N)=1$
we see that the central value $L(G,\chi_{\Delta},1)$ vanishes if $\textrm{sign}(\Delta)\varepsilon_G=-1$, i.e., if $L(G,s)$ has an even functional equation we consider $\Delta<0$ and otherwise $\Delta>0$.

For each of these examples we computed a large set of central derivatives of the twisted $L$-functions with the appropriate $\Delta$ using Sage
and the standard algorithms there which were developed by Dokchitser.
We then fixed a  harmonic weak Maass form with non-zero principal part $P_f$ such that $\xi_{\frac{3}{2}}\left(f \right)$ maps to $G$ under the Shimura lift.
In all cases we took a Poincar\'e series $P_{-\Delta}$ having principal part $q^{-\frac{\Delta}{4N}}$ and computed an initial set of Fourier coefficients for this function using the methods described in the previous section.
We then used the second phase of the algorithm and computed more Fourier coefficients.

Note that for the results in this section, all initial ``phase 1'' computations were all performed using the new Sage package and all further, ``phase 2'', computations were done in Fortran 90.

We would like to give a flavour of the cpu-times involved. The initial computations, using our Sage code, took in all cases approximately 2 hours on a 2.66GHz Xeon processor. On the same processor, the cpu time for a single stretch of phase 2 calculations range between
less than an hour for the smallest discriminant up to several days for the largest discriminant.

As a measure of the accuracy of our computations one can consider the
difference between the coefficients in Tables \ref{tab:11a1-zeros}, \ref{tab:37a1-zeros} and \ref{tab:37b1-zeros}
and the nearest integer (the third column).
To further support the correctness we also list, in Tables \ref{tab:11a1-cminus}, \ref{tab:37a1-cminus} and \ref{tab:37b1-cminus},  normalized coefficients of the non-holomorphic parts,
i.e. $\sqrt{\left|\Delta\right|}c^{-}(\Delta)/\sqrt{\left|\Delta_0\right|}c^{-}(\Delta_0)$
by some fixed non-zero coefficient of index $\Delta_0$.

\subsubsection{11a1}
Here the unique newform of weight two and level 11 is given by
\[G=\eta(\tau)^2 \eta(11\tau)^2=q-2q^2- q^3 + 2q^4 + q^5 +\cdots \in S^{new}_{2}\left(\Gamma_0(11)\right)
\]
and the corresponding $L$-function $L(G,s)$ has an even functional
equation. Using Sage we computed all values of $L'(G,\chi_{\Delta},1)$
for fundamental discriminants $\Delta<0$ such that
$\left(\frac{\Delta}{11}\right)=1$ and $|\Delta|\le19703$.  This set
consists of $2749$ fundamental discriminants and amongst these we
found $14$ discriminants for which $L'(G,\chi_{\Delta},1)$ vanished up
to the numerical precision (see Table \ref{tab:11a1-zeros}).

As a representative for the harmonic weak Maass form in the space
$H_{1/2,\bar\rho}$ corresponding to $G$, we choose the
Poincar\'e series $P_{-5}$ with the principal part
$q^{-\frac{5}{44}}(\frake_7-\frake_{-7})$.  To compute the Fourier coefficients of $P_{-5}$
we used the method described in the previous section with an initial
$\varepsilon=10^{-40}$ and $Y=0.5$, which gave us a truncation point
of $M_0=42$, corresponding to $\Delta$ between $-1847$ and $1885$.
For a short selection of computed values of $c^+(\Delta)$ see Table
\ref{tab:11a1} and for a table of coefficients corresponding to all
vanishing $L'(G,\chi_{\Delta},1)$ see Table \ref{tab:11a1-zeros}.  The
first few normalized ``negative'' coefficients are displayed in Table
\ref{tab:11a1-cminus}. These values should be compared to the list in
\cite[p.\,505]{Sk2}.

\subsubsection{37a1}
Consider the newform of weight two and level 37 which has an odd functional equation. The $q$-expansion is given by
\[
G=q - 2q^2 - 3q^3 + 2q^4 - 2q^5 + 6q^6 - q^7 + 6q^9 + 4q^{10} - 5q^{11} + \cdots \in S^{new}_{2}\left(\Gamma_0(37)\right).
\]
Using Sage we computed all values of $L'(G,\chi_{\Delta},1)$ for fundamental discriminants $\Delta>0$ such that $\left(\frac{\Delta}{37}\right)=1$ and $|\Delta|\le 15000$.
This set consists of $2217$ fundamental discriminants and amongst these we found $8$ discriminants for which $L'(G,\chi_{\Delta},1)$ vanished up to the numerical precision (see Table \ref{tab:37a1-zeros}).
For the corresponding harmonic weak Maass form in $H_{1/2,\rho}$ we took $P_{-3}$, which has a principal part $q^{-\frac{3}{148}}(\frake_{21}+\frake_{21})$.
The initial computation was done in Sage, using $\varepsilon=1\cdot10^{-35}$, which gave a value of $M_0=30$, corresponding to discriminants in the range $-4440 \le \Delta \le 4585$.
For examples of the coefficients $c^{+}(\Delta)$ see Tables  \ref{tab:37a1} and \ref{tab:37a1-zeros}.
The first few normalized ``negative'' coefficients are displayed in Table \ref{tab:37a1-cminus}.

\subsubsection{37b1}
In this case we consider the newform of weight two and level 37 which has an even functional equation. The $q$-expansion is given by
\[
G=q + q^3 - 2q^4 - q^7 - 2q^9 + 3q^{11} + \cdots \in S^{new}_{2}\left(\Gamma_0(37)\right).
\]
Using Sage we computed all values of $L'(G,\chi_{\Delta},1)$ for
fundamental discriminants $\Delta<0$ such that
$\left(\frac{\Delta}{37}\right)=1$ and $|\Delta|\le 12000$.  This set
consists of $1631$ fundamental discriminants and amongst these we
found $15$ discriminants for which $L'(G,\chi_{\Delta},1)$ vanished up
to the numerical precision (see Table \ref{tab:37b1-zeros}).  For the
corresponding harmonic weak Maass form in $H_{1/2,\bar\rho}$ we took
$P_{-12}$, which has a principal part
$q^{-\frac{12}{148}}(\frake_{30}-\frake_{30})$.
The initial computation was
done in Sage, using $\varepsilon=1\cdot10^{-30}$, which gave a value
of $M_0=33$, corresponding to discriminants in the range $-4883 \le
\Delta \le 5029$.  For examples of the coefficients $c^{+}(\Delta)$
see Tables \ref{tab:37b1} and \ref{tab:37b1-zeros}.  The first few
normalized ``negative'' coefficients are displayed in Table
\ref{tab:37b1-cminus}.



\subsection{Conclusions of the numerical experiments for weight two}

In each of the examples of weight two newforms that we studied we saw
agreement with the theorem, i.e.  the coefficients $c^{+}(\Delta)$
(for fundamental discriminants with the appropriate property) were
only algebraic when the corresponding central derivative
$L'(G,\chi_{\Delta},1)$ vanished.  Furthermore, we observed that in
the cases we considered, the algebraic coefficients $c^{+}(\Delta)$
were in fact even rational {\em integers}.

\subsection{Further computations}

To investigate whether a result analogous to Theorem \ref{Lvalues}
also holds for newforms of weight $4$, we computed
$L'(2,G,\chi_{\Delta})$ for all newforms $G$ of weight $4$ on
$\Gamma_0(N)$ whith $5 \le N\le 150$ and fundamental discriminants
$\Delta$ with $|\Delta|\le 300$ and the property that the twisted
$L$-function $L(s,G,\chi_{\Delta})$ has an odd functional equation.
For $5\le N \le 10$ we additionally computed these values for
fundamental discriminants $\Delta$ with $|\Delta|\le5000$.  Amongst
all these values we did not find a single example of a vanishing
derivative.  Even though we did not get any positive case where we
could test the theorem we still wanted to make sure that there was no
easily accesible counter example.

We therefore computed the Fourier coefficients, up to 40 digits
precision, of the associated weight $-\frac{1}{2}$ harmonic Maass form
corresponding to all weight $4$ newforms defined over $\mathbb{Q}$ for
$N$ up to $100$. To test the accuracy (and making sure that the
implementation was correct) we did not only rely on the provable error
bounds, but also checked algebraicity of certain coefficients
corresponding to non-fundamental discriminants. These coefficients
were indeed all found to be integers or rational with fairly small
denominators.  In contrast to this, the Fourier coefficients
corresponding to fundamental discriminants were found not to be similarly ``simple''
rational numbers.

The $L$-value computations were performed in Sage \cite{SA}, using the
included version of Rubinstein's lcalc library \cite{L}.


\subsection{Tables}

\begin{sidewaystable}
\caption[]{$E=11a1$, $P_{-5}\in H_{\frac{1}{2},\bar\rho}$}
\label{tab:11a1}
\begin{tabular}{rll}	
\multicolumn{1}{c}{$\Delta$} &
\multicolumn{1}{c}{$c^{+}(\Delta)$} &
\multicolumn{1}{c}{$L'(G,\chi_{\Delta},1)$}\\
\hline\noalign{\smallskip}
$-7 $&$ \hphantom{-}2.8463370190285980186651576519711751393948073861988\cdot 10^{00} $&$ 1.22556687406888\cdot 10^{00}$\\
$-8 $&$ \hphantom{-}2.5138482002575729165892711124435774460030012341762\cdot 10^{00} $&$ 1.88791354720393\cdot 10^{00}$\\
$-19$&$ \hphantom{-}4.8192428148963255861870924437043119000673519519042\cdot 10^{00} $&$ 7.51391655667451\cdot 10^{00}$\\
$-24$&$ \hphantom{-}5.1018494088339703187507177747597327598232414476402\cdot 10^{00} $&$ 4.02744559024300\cdot 10^{00}$\\
$-35$&$ \hphantom{-}4.7515892636101723769649079639675162004017245362399\cdot 10^{00} $&$ 7.64786334637073\cdot 10^{00}$\\
$-39$&$            -1.6466690697010481166272091028219327677356442914804\cdot 10^{01} $&$ 2.97721567216550\cdot 10^{00}$\\
$-40$&$ \hphantom{-}1.1470941388138074683747768314723689292860539962900\cdot 10^{01} $&$ 5.58789208952436\cdot 10^{00}$\\
$-43$&$-1.7622439638503327722737780360046423237568367805048\cdot 10^{01} $&$ 1.18814465355690\cdot 10^{01}$\\
$-51$&$ \hphantom{-}2.0736222999878741718629718432682995552582880786065\cdot 10^{01} $&$ 1.30416363302768\cdot 10^{01}$\\
$-52$&$ \hphantom{-}1.5723528683914990387103216700146317562411438497615\cdot 10^{01} $&$ 5.14853759817659\cdot 10^{00}$\\
$-68$&$ \hphantom{-}9.6889673322938493992006043404127469979247370067926\cdot 10^{00} $&$ 3.80344864881298\cdot 10^{00}$\\
$-79$&$ \hphantom{-}1.7557351755436160739388564340027760291317089229254\cdot 10^{01} $&$ 4.75620653690677\cdot 10^{00}$\\
$-83$&$            -7.1767664383427675609861242907417950544611683162859\cdot 10^{01} $&$ 6.43843846621214\cdot 10^{00}$\\
$-84$&$ \hphantom{-}6.1666200626587315159968126799603650525586539365601\cdot 10^{01} $&$ 6.53746327159376\cdot 10^{00}$\\
$-87$&$ -7.7230036424433334541697484050338439023979647483280\cdot 10^{01} $&$ 2.35584785481347\cdot 10^{00}$\\
$-95$&$ \hphantom{-}7.8467572084064151556661839046504144426199227897994\cdot 10^{01} $&$3.03660486030085\cdot 10^{00}$\\
$-811$&$\hphantom{-}3.0046247983067285336553431175489765847382042907105\cdot 10^{06} $&$1.25949136911120\cdot 10^{01}$\\
$-820$&$-6.0493754250387304262091147332158578046510749019315\cdot 10^{06} $&$1.19119437485937\cdot 10^{01}$\\
$-824$&$-5.7985199999999999999999999999999999999999999999999\cdot 10^{06} $&$-6.6\cdot 10^{-24}$\\
$-827$&$ \hphantom{-}1.8535489407871222859528059423736067736521222528554\cdot 10^{06} $&$1.60961273159300\cdot 10^{01}$\\
$-831$&$-6.7911392225835416131083026699411310608420151994986\cdot 10^{06} $&$3.36744068632019\cdot 10^{00}$\\
$-996$&$-3.5516294505685820400211045047063422129168082892941\cdot 10^{07}$ & $1.15828152096335\cdot 10^{01}$\\
$-1003$&$ \hphantom{-}1.0811934742079303073802766406181476437608668928409\cdot 10^{07} $&$2.53076681967579\cdot 10^{01}$\\
$-1007$&$-3.9469248000000000000000000000000000000000000000000\cdot 10^{07} $& $-1.1\cdot 10^{-22}$\\
$-1011$&$-3.7685140824429636488934775010060106547341275101054\cdot 10^{07} $&$ 1.84592490209627\cdot 10^{01}$\\
$-1019$&$ \hphantom{-}3.3790315957549749442218769650593817997888628818338\cdot 10^{07}$ &$ 1.68145450009782\cdot 10^{01}$ \\

\end{tabular}
\end{sidewaystable}

\begin{table}
\caption[]{$E=11a1$, $P_{-5} \in H_{\frac{1}{2},\bar\rho}$}
\label{tab:11a1-zeros}
\begin{tabular}{llc}	
\multicolumn{1}{c}{$\Delta$} &
\multicolumn{1}{c}{$c^{+}(\Delta)$} &
\multicolumn{1}{c}{$\left|c^{+}(\Delta)-[c^{+}(\Delta)]\right|$} \\
\hline\noalign{\smallskip}
$-824 $&$ -5798520  $&$ 3.0\cdot 10^{-76}$   \\
$-1799 $&$-2708450784 $&$ 2.7\cdot 10^{-46}$\\
$-4399 $&$-68135748249936640 $&$ 2.3\cdot 10^{-21}$ \\
$-8483 $&$ \hphantom{-}214445760716391388216704 $&$ 9.1\cdot 10^{-28}$  \\
$-11567 $&$ -12412267149099919205092899456 $&$ 1.6 \cdot 10^{-25} $\\
$-14791 $&$ \hphantom{-} 66850179291021019012709832099520 $&$ 3.1\cdot 10^{-30}$\\
$-15487  $&$ -478732239405182448762415030881280 $&$ 5.6\cdot 10^{-30}$ \\
$-15659  $&$ -804489814454597618648064770159415 $&$ 6.8\cdot 10^{-30}$ \\
$-15839  $&$ -1162122495004344641799524116135680 $&$ 7.0\cdot 10^{-30}$ \\
$-16463  $&$  \hphantom{-}4542575922533728228643934862230144 $&$ 1.5\cdot 10^{-30} $\\
$-17023  $&$ -23302350713109514450879400185948800 $&$ 2.0\cdot 10^{-29}$\\
$-17927  $&$  \hphantom{-}110133238181959291703634158808374784 $&$ 1.2\cdot 10^{-29}$\\
$-18543  $&$  \hphantom{-}464726791864282489334104058164482624 $&$  1.9\cdot 10^{-29}$\\
\end{tabular}
\end{table}

\begin{table}
\caption[]{$E=11a1$, $P_{-5}\in H_{\frac{1}{2},\bar\rho}$. Coefficients are scaled by
$ c^{-}(1)$.}
\label{tab:11a1-cminus}
\begin{tabular}{llc}	
\multicolumn{1}{c}{$\Delta$} &
\multicolumn{1}{c}{$\sqrt{\Delta}\,c^{-}(\Delta)$} &
\multicolumn{1}{c}{$\left|c^{-}(\Delta)-[c^{-}(\Delta)]\right|$} \\
\hline\noalign{\smallskip}
$ 4 $&$-3$  &$2.0 \cdot 10^{-100}$\\
$ 5 $&$\hphantom{-}5$  &$2.1 \cdot 10^{-99\hphantom{0}}$\\
$ 9 $&$-2$ & $1.7 \cdot 10^{-100}$\\
$ 12$&$\hphantom{-}5$  &$8.0 \cdot 10^{-100}$\\
$ 16$&$\hphantom{-}4$  &$1.5 \cdot 10^{-99\hphantom{0}}$\\
$ 20$&$\hphantom{-}5$  &$1.1 \cdot 10^{-100}$\\
$ 25$&$\hphantom{-}0$  &$1.0 \cdot 10^{-100}$\\
$ 36$&$\hphantom{-}6$  &$1.0 \cdot 10^{-99\hphantom{0}}$\\
$ 37$&$\hphantom{-}5$  &$4.2 \cdot 10^{-99\hphantom{0}}$\\
$ 45$&$\hphantom{-}0$  &$6.4 \cdot 10^{-99\hphantom{0}}$\\
\end{tabular}
\end{table}

\begin{sidewaystable}
\caption[]{$E=37a1$, $P_{-3}\in H_{\frac{1}{2},\rho}$}
\label{tab:37a1}
\begin{tabular}{rll}
\multicolumn{1}{c}{$\Delta$} &
\multicolumn{1}{c}{$c^{+}(\Delta)$} &
\multicolumn{1}{c}{$L'(G,\chi_{\Delta},1)$}\\
\hline\noalign{\smallskip}
$ 1  $&$-2.8176178498959956879756075537515493438922975370716\cdot 10^{-01}$ & $3.05999773834052\cdot 10^{-01}$\\
$ 12 $&$-4.8852723826201225228227029607337071669095284814788\cdot 10^{-01}$ & $4.29861479867736\cdot 10^{00}$\\
$ 21 $&$-1.7273925723265275652082007397068992218426924398791\cdot 10^{-01}$ & $9.00238680032537\cdot 10^{00}$\\
$ 28 $&$\hphantom{-}6.7819399530394779828450578400669420938246859928076\cdot 10^{-01}$ & $4.32726024966011\cdot 10^{00}$\\
$ 33 $&$\hphantom{-}5.6630232015906998168220545669245622604190884430064\cdot 10^{-01}$ & $3.62195679113882\cdot 10^{00}$\\
$ 37 $&$-9.1326561374611652958506448407204050631184401026129\cdot 10^{-01}$ & $3.47328771649229 \cdot 10^{00}$\\
$ 40 $&$\hphantom{-}4.0098509269543637915254766073122850557290259963615\cdot 10^{-01}$ & $3.70588717878444\cdot 10^{00}$\\
$ 41 $&$\hphantom{-}6.5637495744757231959699415722023547525400239084778\cdot 10^{-01}$ & $5.93680171871573\cdot 10^{00}$\\
$ 44 $&$\hphantom{-}9.6886404434506397321859573425794267139455920322171\cdot 10^{-01}$ & $1.01334656625280\cdot 10^{01}$\\
$ 53 $&$-5.6688852568232517859984506645723944339238503996096\cdot 10^{-01}$ & $2.61746665637296\cdot 10^{01}$\\
$ 65 $&$-6.0328072889521477971996798071175156059595671497733\cdot 10^{-01}$ & $7.67818286326206\cdot 10^{00}$\\
$ 73 $&$\hphantom{-}3.4874711835362408853804154923777452565842552803223\cdot 10^{-01}$ &$2.92507284795068\cdot 10^{00}$\\
$ 77 $&$\hphantom{-}2.2699132373705254600799448087564660809534768699467\cdot 10^{-01}$ &$3.42067600398534\cdot 10^{10}$\\
$ 85 $&$-7.6894617048676272061865441758881289699552927122551\cdot 10^{-01}$ & $9.90133670369251\cdot 10^{00}$\\
$ 1481$&$-3.2715595098273932057423414526408419506801164996884\cdot 10^{00}$ &$5.26994449124823\cdot 10^{00}$ \\
$ 1484$&$-1.3432792297590353562651264178555980321660674399890\cdot 10^{01}$ &$3.86746474997364\cdot 10^{01}$\\
$ 1489$& $\hphantom{-}8.9999999999999999999999999999999999999999999999999\cdot 10^{00}$ &$-3.7\cdot 10^{-23}$\\
$ 1496$&$\hphantom{-}1.1199440423162819213593329208218112792285448658029\cdot 10^{01}$ &$2.27616829409607\cdot 10^{01}$\\
$ 1501$&$-5.8188238119388864901078937905792951783427273426771\cdot 10^{02}$ &$ 6.06007663972706\cdot 10^{00}$\\
$ 4376$&$-3.6731327299348159991042234350611468700145535059868\cdot 10^{02}$ &$2.03155740209437\cdot 10^{01}$\\
$ 4377$&$-5.0062522276143084997015960658866819832113068397294\cdot 10^{02}$ &$2.27150950159608\cdot 10^{00}$\\
$ 4393$&$\hphantom{-}6.6000000000000000000000000000000000000000000001468\cdot 10^{01}$ &$5.8\cdot 10^{-23}$\\
$ 4396$&$-2.3023069110811173762943326075771836710063196221488\cdot 10^{02}$ &$2.00437958330233\cdot 10^{00}$\\
$ 4412$&$-3.1500483730098996665306117169085504925545562420809\cdot 10^{02}$ &$3.73011222569745\cdot 10^{01}$\\
\end{tabular}
\end{sidewaystable}

\begin{table}
\caption[]{$E=37a1$, $P_{-3}\in H_{\frac{1}{2},\rho}$}
\label{tab:37a1-zeros}
\begin{tabular}{llc}	
\multicolumn{1}{c}{$\Delta$} &
\multicolumn{1}{c}{$c^{+}(\Delta)$} &
\multicolumn{1}{c}{$\left|c^{+}(\Delta)-[c^{+}(\Delta)]\right|$} \\
\hline\noalign{\smallskip}
$1489$&$\hphantom{-}9 $    &$ 1.6\cdot 10^{-72}$\\
$4393$&$\hphantom{-}66 $   &$  1.5\cdot 10^{-45}$\\
$5116$&$-746$  &$ 8.5\cdot 10^{-23}$\\
$5281$&$\hphantom{-}153$  &$ 8.2\cdot 10^{-23}$\\
$5560$&$-1124$ &$  1.2\cdot 10^{-22}$\\
$5761$&$-974$  &$ 1.1\cdot 10^{-22}$\\
$6040$&$-1404$ &$  4.2\cdot 10^{-23}$\\
$6169$&$\hphantom{-}336$  &$  1.1\cdot 10^{-22}$\\
\end{tabular}
\end{table}

\begin{table}
\caption[]{$E=37a1$, $P_{-3}\in H_{\frac{1}{2},\rho}$. Coefficients are scaled by $\sqrt{3}\,c^{-}(-3)$.}
\label{tab:37a1-cminus}
\begin{tabular}{ccc}	
\multicolumn{1}{c}{$\Delta$} &
\multicolumn{1}{c}{$\sqrt{|\Delta|}\,c^{-}(\Delta)$} &
\multicolumn{1}{c}{$\left|c^{-}(\Delta)-[c^{-}(\Delta)]\right|$} \\
\hline\noalign{\smallskip}
$ -4$&$\hphantom{-}1$ &$4.0 \cdot 10^{-84}$\\
$ -7$&$-1$ &$5.0 \cdot 10^{-84}$\\
$-11$&$\hphantom{-}1$ &$4.5 \cdot 10^{-84}$\\
$-12$&$-1$ &$2.0 \cdot 10^{-84}$\\
$-16$&$-2$ &$1.1 \cdot 10^{-83}$\\
$-27$&$-3$ &$1.3 \cdot 10^{-83}$\\
$-28$&$\hphantom{-}3$ &$1.4 \cdot 10^{-83}$\\
$-36$&$-2$ &$1.0 \cdot 10^{-83}$\\
$-40$&$\hphantom{-}2$ &$3.6 \cdot 10^{-85}$\\
$-44$&$-1$ &$1.1 \cdot 10^{-83}$\\
\end{tabular}
\end{table}

\begin{sidewaystable}
\caption[]{$E=37b1$, $P_{-12}\in H_{\frac{1}{2},\bar\rho}$}
\label{tab:37b1}
\begin{tabular}{rll}
\multicolumn{1}{c}{$\Delta$} &
\multicolumn{1}{c}{$c^{+}(\Delta)$} &
\multicolumn{1}{c}{$L'(G,\chi_{\Delta},1)$}\\
\hline\noalign{\smallskip}
$ -3 $&$\hphantom{-}1.0267149116920354474451980654263083626994977508118 \cdot 10^{00}$&$1.47929949208 \cdot 10^{00}$\\
$ -4 $&$\hphantom{-}1.2205364009670316625279102409757685190938519711297 \cdot 10^{00}$&$1.81299789722 \cdot 10^{00}$\\
$ -7 $&$\hphantom{-}1.6900297463200076214148752932012403965170838158011 \cdot 10^{00}$&$2.11071898018 \cdot 10^{00}$\\
$-11 $&$\hphantom{-}5.8849982354849175483779961900586424239744874288522 \cdot 10^{-01}$&$3.65679089534 \cdot 10^{00}$\\
$-40 $&$\hphantom{-}1.2669706585839831188366862729215921230412462308976 \cdot 10^{00}$&$4.16362898338 \cdot 10^{00}$\\
$-47 $&$\hphantom{-}3.0756790552662277517712909874001702657447701023621 \cdot 10^{00}$&$5.26739088546 \cdot 10^{00}$\\
$-67 $&$\hphantom{-}2.1608356105538234382282266707128748893591830597455 \cdot 10^{00}$&$4.98143961845 \cdot 10^{00}$\\
$-71 $&$-1.5945418432752378367351454423028372659103804842831\cdot 10^{00}$&$5.33295381308 \cdot 10^{00}$\\
$-83 $&$\hphantom{-}2.9631171578917930530100644900469583789690213329975 \cdot 10^{00}$&$7.30522465208 \cdot 10^{00}$\\
$-84 $&$-3.8773494709413749500399075799371212202544017987791\cdot 10^{00}$&$1.00026475317 \cdot 10^{01}$\\
$-95 $&$-2.6554862688645143792016861758519887392731139540185 \cdot 10^{00}$&$5.83606039003 \cdot 10^{00}$\\
$-132$&$\hphantom{-}4.1944733541115532186541330550136737538249181082859 \cdot 10^{00}$&$9.99216716471 \cdot 10^{00}$\\
$-136$&$-4.8392675993443437829864850814885823635034770966657 \cdot 10^{00}$&$5.73824076491 \cdot 10^{00}$\\
$-139$&$-5.9999999999999999999999999999999999999999999999991 \cdot 10^{00}$&$-8.5 \cdot 10^{-23}$\\
$-151$&$-8.3135688179267692046624844818371994826339638923811 \cdot 10^{-01}$&$6.69750855159 \cdot 10^{00}$\\
$-152$&$\hphantom{-}4.3274351625459058613812696410805017025617476195953 \cdot 10^{00}$&$7.95190347996 \cdot 10^{00}$\\
$-811$&$-1.4731293182498551151700589944493338505308298027148 \cdot 10^{02}$&$5.32436617837 \cdot 10^{00}$\\
$-815$&$\hphantom{-}1.2194410312093092058885476868302234805268383148338 \cdot 10^{02}$&$4.74925836935 \cdot 10^{00}$\\
$-823$&$\hphantom{-}3.1200000000000000000000000000000000000000000000000 \cdot 10^{02}$&$-1.5 \cdot 10^{-23}$\\
$-824$&$-3.2299860660409750567356931348586086493552010570382 \cdot 10^{02}$&$1.75028741141 \cdot 10^{01}$\\
$-835$&$-2.4035736526655124690110045874885626910322384359422 \cdot 10^{02}$&$8.64359690730 \cdot 10^{00}$\\
\end{tabular}
\end{sidewaystable}

\begin{table}
\caption[]{$E=37b1$, $P_{-12}\in H_{\frac{1}{2},\bar\rho}$}
\label{tab:37b1-zeros}
\begin{tabular}{llc}	
\multicolumn{1}{c}{$\Delta$} &
\multicolumn{1}{c}{$c^{+}(\Delta)$} &
\multicolumn{1}{c}{$\left|c^{+}(\Delta)-[c^{+}(\Delta)]\right|$} \\
\hline\noalign{\smallskip}
$ -139 $  &  $ -6 $  &  $ 1.5 \cdot 10^{-85} $ \\
 $ -823 $  &  $\hphantom{-}312 $  &  $ 9.1 \cdot 10^{-80} $ \\
 $ -2051 $  &  $ -26724 $  &  $ 1.0 \cdot 10^{-67} $ \\
 $ -2599 $  &  $\hphantom{-}122048 $  &  $ 3.4 \cdot 10^{-63} $ \\
 $ -3223 $  &  $ -472416 $  &  $ 3.2 \cdot 10^{-57} $ \\
 $ -3371 $  &  $ -674712 $  &  $ 7.4 \cdot 10^{-56} $ \\
 $ -5227 $  &  $\hphantom{-}5816 $  &  $ 5.5 \cdot 10^{-31} $ \\
 $ -5307 $  &  $ -5192 $  &  $ 4.6 \cdot 10^{-31} $ \\
 $ -6583 $  &  $ -13320 $  &  $ 4.6 \cdot 10^{-31} $ \\
 $ -7892 $  &  $ -79552 $  &  $ 1.2 \cdot 10^{-30} $ \\
 $ -7951 $  &  $\hphantom{-}28152 $  &  $ 4.0 \cdot 10^{-31} $ \\
 $ -9112 $  &  $ -224548 $  &  $ 1.6 \cdot 10^{-30} $ \\
 $ -9715 $  &  $\hphantom{-}236934 $  &  $ 2.8 \cdot 10^{-32} $ \\
 $ -11444 $  &  $ -1437956 $  &  $ 2.0 \cdot 10^{-33} $ \\
 $ -11651 $  &  $\hphantom{-}563716 $  &  $ 7.0 \cdot 10^{-34} $ \\
\end{tabular}
\end{table}

\begin{table}
\label{tab:37b1-cminus}
\caption[]{$E=37b1$, $P_{-12}\in H_{\frac{1}{2},\bar\rho}$. Coefficients are scaled by $c^{-}(1)$.}
\begin{tabular}{ccc}
\multicolumn{1}{c}{$\Delta$} &
\multicolumn{1}{c}{$\sqrt{\Delta}\,c^{-}(\Delta)$} &
\multicolumn{1}{c}{$\left|c^{-}(\Delta)-[c^{-}(\Delta)]\right|$} \\
\hline\noalign{\smallskip}
$ 4 $&$-1$ &$1.6 \cdot 10^{-85}$\\
$ 9 $&$\hphantom{-}0$ &$3.2 \cdot 10^{-85}$\\
$ 12$&$\hphantom{-}3$ &$9.6 \cdot 10^{-85}$\\
$ 16$&$-2$ &$1.6 \cdot 10^{-85}$\\
$ 21$&$\hphantom{-}3$ &$2.7 \cdot 10^{-85}$\\
$ 25$&$-1$ &$1.7 \cdot 10^{-85}$\\
$ 28$&$\hphantom{-}3$ &$3.9 \cdot 10^{-85}$\\
$ 33$&$\hphantom{-}3$ &$3.6 \cdot 10^{-85}$\\
$ 36$&$\hphantom{-}0$ &$3.9 \cdot 10^{-85}$\\
$ 40$&$\hphantom{-}0$ &$5.6 \cdot 10^{-85}$\\
\end{tabular}
\end{table} 

\end{document}